\documentclass[a4paper,12pt]{article}
\usepackage{graphicx} 
\usepackage[utf8]{inputenc}
\usepackage[T2A]{fontenc}
\usepackage{amsfonts}
\usepackage{amssymb,amsmath,amsthm}
\usepackage{bm}
\usepackage{tikz}
\usepackage{tikz-cd}
\usepackage{amsthm}
\usepackage{import}
\usepackage[left=2cm,right=2cm,top=2cm,bottom=2cm,bindingoffset=0cm]{geometry}

\usepackage{amscd}
\usepackage[all]{xy}

\def\ve{\varepsilon}
\newtheorem{theorem}{Theorem}
\newtheorem*{definition*}{Definition}

\newtheorem{lemma}{Lemma}
\newtheorem*{claim*}{Remark}
\newtheorem{proposition}{Proposition}

\def\beq{\begin{equation}}
	\def\eeq{\end{equation}}
\def\beqq{\begin{equation*}}
	\def\eeqq{\end{equation*}}
\def\bg{{\bm{\gamma}}}
\def\bbg{\widehat{\bg}}

\def\bp{\mathbf{p}}	
\def\bx{\mathbf{x}}	
\def\d{\partial}
\def\g{\gamma}
\def\gam{\nu}
\def\bgam{\bm{\nu}}

\def\C{\mathbb{C}}
\def\R{\mathbb{R}}
\def\tt{\theta}

\newcommand{\rf}[1]{(\ref{#1})}

\def\cbk{\color{black}}

\def\ii{\operatorname{{\mathbf i}}}
\def\jj{\operatorname{{\mathbf j}}}
\def\ss{\operatorname{{\mathbf s}}}
\begin{document}
	
	\begin{center} {\bf\large  Zhelobenko - Stern 
			formulas and  $B_n$ Toda wave functions}\bigskip
		
		{\bf A.Galiullin$^{\circ}$,  S. Khoroshkin$^{\circ\ast}$, M.Lyachko$^{\circ}$
		}\medskip\\
		$^\circ${\it National Research University Higher School of Economics, Myasnitskaya 20, \\Moscow, 101000, Russia;}\smallskip\\
		$^\ast${\it 
		Institute for Information Transmission Problems RAS (Kharkevich Institute), \\Bolshoy Karetny per. 19, Moscow, 127994, Russia}
	
\end{center}

\begin{abstract} Using Zhelobenko - Stern formulas for the action of the generators of orthogonal Lie algebra in corresponding Gelfand - Tsetlin basis, we derive Mellin - Barnes presentations for the wave functions of $B_n$ Toda lattice. They are in accordance with Iorgov - Shadura formulas.
	\end{abstract}

	\section{Introduction} In the paper \cite{GKL} Gerasimov, Kharchev and Lebedev applied the famous  formulas \cite{GT} for the action of generators of general Lie algebra $gl(n)$ in Gelfand - Tsetlin basis of irreducible finite - dimensional representations of general linear group $GL(n,\C)$ 
	 to obtain Mellin - Barnes presentation of the wave functions of open $A_n$ Toda chain. Using Gelfand -Tsetlin formulas, they constructed an infinite - dimensional representation of Lie algebra $gl(n,\C)$ in the space of meromorphic functions on $n(n-1)/2$ variables, found there two dual Whittaker vectors and realized, according Kostant theory \cite{K}, the Toda wave function as certain  matrix element in this representation.
	  The same formulas were earlier established by Kharchev and Lebedev in the technique of Yang-Baxter formalism \cite{KL}. 
	  
	  Besides, Toda wave functions admits another presentation, known by the name Gauss - Givental  by means of integrals over spatial variables. It was found first in \cite{Giv}. Gauss - Givental presentation   was then derived in \cite{GLO} for  wave functions of Toda systems related to $B_n$, $C_n$, $D_n$ root systems.
	    
      Mellin transform of these formulas was computed in \cite{KK}. However, the formulas presented in \cite{KK} differ from that of \cite{GKL} and are not satisfactory by several reasons. In particular, one cannot find in this presentation  Sklyanin measure and thus the results  of \cite{KK} cannot be used to establish the completeness and orthogonality of the wave function and develop the corresponding integral transform.
      
      The goal of this paper is to try to fill this gap at least partially for $B_n$ Toda system using representation theoretical tools similar to that of \cite{GKL}. The only known result in this direction is the paper \cite{IS} of Iorgov and Shadura where they constructed $B_n$ wave function by its decomposition over related $A_n$ Toda wave function. As well as in \cite{KL}, this work was done in a framework of the Yang - Baxter formalism.
      
      Our starting point is an analog of Gelfand - Tsetlin formulas for orthogonal groups published without a proof by Zhelobenko and Stern in \cite{ZS}. These formulas look much more complicated compared to \cite{GT} and we did not find numerous applications of them in the literature. However, after their check we constructed infinite dimensional representation of the orthogonal Lie algebra and found there two dual Whittaker vectors. With their help we constructed the integrals, presenting $B_n$ wave functions in which we see all expected ingredients of Sklyanin measure. The resulting formula can be presented as an iterative procedure in two ways.
      
       Firstly, it is an iterative procedure over the rank of orthogonal group and this is probably the most interesting result of this paper. Each step can be interpreted as an action of the raising integral operator, which kernel is itself an integral over intermediate variables. Such type of structure we also observe in Gauss - Givental representation \cite[(1.74)]{GLO}.
       Second, we can consider two successive iterative integrals combining them in other parity. Then intermediate step becomes precisely a degeneration of $B_n$ Gustafson integral and can be explicitly evaluated. In this way we arrive at Iorgov - Shadura formula.
       
       Note  two subtle points of our construction. First, Zhelobenko-Stern formulas are written for the generators of orthogonal Lie algebras in their orthogonal realization, while  Whittakker vectors refer to simple root generators. An existence of Whittakker vectors in a factorized form was not evident from the beginning. By the same reasons  action of the Cartan subalgebra in corresponding infinite dimensional representation cannot be written,  contrary to $gl(n)$, in terms of multiplications by linear functions. Fortunately, it is so for the action on Whittaker vectors.
       
         Despite the fact that Zhelobenko -Stern formulas are written uniformly for all orthogonal Lie algebras, we managed to find Whittaker vectors only for Lie algebras $so(2n-1)$. More precisely, we managed firstly to found factorized expressions for 'degenerate' Whittaker vectors for $so(2n)$, for which one of the simple generators act by zero, so that they are essentially Whittaker vectors for embedded $gl(n)$ Lie algebra. Restricting these vectors to $so(2n-1)$ we get 'nondegenerate' Whittaker vectors for this Lie algebra.

	\section{Gelfand-Tsetlin type representation}
	\subsection{Zhelobenko-Stern formulas}
	It is well known that each irreducible representation of the orthogonal groups $SO(2n+1)$ and $SO(2n)$ is parametrized by its signature, given by  ordered sequences of integers, respectively
	\beq\label{0a} \begin{split} p_1\geq p_2\geq\ldots p_{n-1}\geq p_n\geq 0\\
		p_1\geq p_2\geq\ldots p_{n-1}\geq | p_n|\end{split}\eeq
	and the restriction of irreducible representation of $SO(2n+1)$ to $SO(2n)$ has simple spectrum described by all signatures $q_1,\ldots q_n$, satisfying interleaving inequalities
	\beq \label{0b}p_1\geq q_1\geq p_2\geq\ldots q_{n-1}\geq p_n\geq q_n\geq - p_n.\eeq
	Analogously,  the restriction of irreducible representation of $SO(2n)$ to $SO(2n-1)$ has simple spectrum described by all signatures $q_1,\ldots q_n$, satisfying interleaving inequalities
	\beq\label{0c} p_1\geq q_1\geq p_2\geq\ldots q_{n-1}\geq p_n.\eeq
	This enables one to construct an orthogonal basis of irreducible representation of the orthogonal group $SO(n)$ parametrised by Gelfand-Tsetlin tableaux  
	\begin{equation} \label{1}
		\mathbf{p} = \left(
		\begin{array}{cccc}
			p_{n-1,1} & p_{n-1,2} & \ldots & p_{n-1,[\frac{n}{2}]}\\
			\vdots & \vdots & \rotatebox{45}{\ldots} \\
			p_{3,1} & p_{3,2} \\
			p_{2,1} \\
			p_{1,1} 
		\end{array}
		\right) = \left(
		\begin{array}{c}
			\mathbf{p}_{n-1} \\
			\vdots \\
			\mathbf{p}_3 \\
			\mathbf{p}_2 \\
			\mathbf{p}_{1} 
		\end{array}
		\right)
	\end{equation} The upper row  $\bp_{n-1}$
	indicates the signature of the irreducible representation of $SO(n)$ and is fixed for all tableaux  parametrising its basic vectors, the second row indicates the signature of the restriction to $SO(n-1)$, etc., the integer $p_{11}$ indicates the irreducible of $SO(2)$. 
	Integers $p_{ij}$ should satisfy the row by row interleaving inequalities \rf{0a}--\rf{0b}.
	
	It is natural to shift the signatures by the corresponding half sum of positive roots of the related root system, that is, we set
	\beq\label{1a} \begin{split} m_{2n,k}&= p_{2n,k}+ (n-k)+\frac{1}{2},\\
		m_{2n-1,k}&=p_{2n-1,k}+ (n-k)\end{split}\eeq
	Zhelobenko and Stern \cite[Chapter II, Section 5.8]{ZS} presented without a proof a precise expression for the matrix elements of generators of the Lie algebra $so(n)$ in the corresponding orthogonal basis. 
	
	The Lie algebra $so(n)$ is generated, as a vector space, by elements 
	$$I_{ij}=e_{ij}-e_{ji},\qquad i>j $$
	As a Lie algebra, it is generated by elements $$I_{i+1,i}, \qquad i=1,\ldots, n-1$$
	with defining relations  
	\beq \label{2} \begin{split}
		& \left[ I_{i+1,i}, \left[ I_{i+2,i+1}, I_{i+1,i}\right] \right] = I_{i+2,i+1} \quad i=1,\ldots,n-2 \\ &
		\left[ I_{i+2,i+1}, \left[ I_{i+1,i}, I_{i+2,i+1}\right] \right] = I_{i+1,i} \quad i=1,\ldots,n-2 \\ &
		\left[ I_{i+1,i}, I_{j+1,j} \right] = 0 \quad |i-j| > 1
	\end{split} \eeq
	After a renormalization, eliminating square roots in the coefficients and correcting misprints, their formulas look like
	\beq\label{3}\begin{split}
		I_{2k+1,2k} & = -\sum_{j=1}^{k} \frac{\prod\limits_{r=1}^{k-1} (m_{2k-2,r} + m_{2k-1,j} +\frac{1}{2} ) \prod\limits_{r=1}^{k} (m_{2k-1,j} - m_{2k,r} + \frac{1}{2})}{2\prod\limits_{r \neq j} (m_{2k-1,j} - m_{2k-1,r})(m_{2k-1,j} + m_{2k-1,r}+1)} e^{ \partial_{m_{2k-1,j}}} - \\ & - \sum_{j=1}^{k} \frac{\prod\limits_{r=1}^{k-1} (m_{2k-1,j} - m_{2k-2,r}-\frac{1}{2} ) \prod\limits_{r=1}^{k} (m_{2k,r} + m_{2k-1,j} - \frac{1}{2})}{2\prod\limits_{r \neq j} (m_{2k-1,j} - m_{2k-1,r})(m_{2k-1,j} + m_{2k-1,r}-1)} e^{- \partial_{m_{2k-1,j}}}
	\end{split}\eeq
	\beq\label{4}\begin{split}
		I_{2k+2,2k+1} & =  \sum_{j=1}^{k} \frac{\prod\limits_{r=1}^{k+1} ((m_{2k,j}+\frac{1}{2})^2 - m_{2k+1,r}^2)}{2m_{2k,j} (m_{2k,j}+\frac{1}{2}) \prod\limits_{r \neq j} (m^2_{2k,j} - m^2_{2k,r})} e^{ \partial_{m_{2k,j}}} + \\ & +  \sum_{j=1}^{k} \frac{\prod\limits_{r=1}^{k} ((m_{2k,j}-\frac{1}{2})^2 - m_{2k-1,r}^2)}{2m_{2k,j} (m_{2k,j}-\frac{1}{2}) \prod\limits_{r \neq j} (m^2_{2k,j} - m^2_{2k,r})} e^{-\partial_{m_{2k,j}}} + \\ & +i  \frac{\prod\limits_{r=1}^{k} m_{2k-1,r} \prod\limits_{r=1}^{k+1} m_{2k+1,r}}{\prod\limits_{r=1}^{k}(m_{2k,r}+\frac{1}{2}) (m_{2k,r}-\frac{1}{2})}
	\end{split}\eeq
	Here 
	the operators $e^{\pm\d m_{kj}}$ are operators of  shifts of the entries of Gelfand -Tsetlin tableau:
	the operator $e^{\pm\d m_{kj}}$ changes $m_{kj}$ by $m_{kj}\pm 1$ (and respectively  $p_{kj}$ by $p_{kj}\pm 1$).
	We can extend the RHS of  relations \rf{3} and  \rf{4} to arbitrary complex parameters $m_{ij}$ and regard them as operators acting in the space of rational functions on $m_{ij}$.

	\subsection{Representation in meromorphic functions }
	Following \cite{GKL}, we renormalize the variables
	\beq m_{ij}=\frac{\gam_{ij}}{ic} \eeq
	in order to have an additional scaling variable in the representation. Then we have
	\beq\label{5}\begin{split}
		I_{2k+1,2k} & = -\frac{1}{ic}\sum_{j=1}^{k} \frac{\prod\limits_{r=1}^{k-1} (\gam_{2k-2,r} + \gam_{2k-1,j} +\frac{{ic}}{2} ) \prod\limits_{r=1}^{k} (\gam_{2k-1,j} - \gam_{2k,r} + \frac{{ic}}{2})}{2\prod\limits_{r \neq j} (\gam_{2k-1,j} - \gam_{2k-1,r})(\gam_{2k-1,j} + \gam_{2k-1,r}+{ic})} e^{ {ic}\partial_{\gam_{2k-1,j}}} - \\ & -\frac{1}{{ic}} \sum_{j=1}^{k} \frac{\prod\limits_{r=1}^{k-1} (\gam_{2k-1,j} - \gam_{2k-2,r}-\frac{ic}{2} ) \prod\limits_{r=1}^{k} (\gam_{2k,r} + \gam_{2k-1,j} - \frac{{ic}}{2})}{2\prod\limits_{r \neq j} (\gam_{2k-1,j} - \gam_{2k-1,r})(\gam_{2k-1,j} + \gam_{2k-1,r}-{ic})} e^{-{ic} \partial_{\gam_{2k-1,j}}}
	\end{split}\eeq
	\beq\label{6}\begin{split}
		iI_{2k+2,2k+1} & =  \frac{1}{c}\sum_{j=1}^{k} \frac{\prod\limits_{r=1}^{k+1} ((\gam_{2k,j}+\frac{{ic}}{2})^2 - \gam_{2k+1,r}^2)}{2\gam_{2k,j} (\gam_{2k,j}+\frac{{ic}}{2}) \prod\limits_{r \neq j} (\gam^2_{2k,j} - \gam^2_{2k,r})} e^{{ic} \partial_{\gam_{2k,j}}} + \\ & +  \frac{1}{c}\sum_{j=1}^{k} \frac{\prod\limits_{r=1}^{k} ((\gam_{2k,j}-\frac{{ic}}{2})^2 - \gam_{2k-1,r}^2)}{2\gam_{2k,j} (\gam_{2k,j}-\frac{{ic}}{2}) \prod\limits_{r \neq j} (\gam^2_{2k,j} - \gam^2_{2k,r})} e^{-{ic}\partial_{\gam_{2k,j}}} + \\ & -  \frac{1}{ic}\frac{\prod\limits_{r=1}^{k} \gam_{2k-1,r} \prod\limits_{r=1}^{k+1} \gam_{2k+1,r}}{\prod\limits_{r=1}^{k}(\gam_{2k,r}+\frac{{ic}}{2}) (\gam_{2k,r}-\frac{{ic}}{2})}
	\end{split}\eeq 
	\begin{proposition} The operators \rf{3} and \rf{4} satisfy the defining relations \rf{2} of the generators of orthogonal Lie algebra $so(n)$  		
	\end{proposition} 
	Surely, this statement follows from its validity in finite dimensional representations, since the relations are then satisfied on sufficiently many integer points. However, since the proof of the formulas is missing in \cite{ZS}, we checked the defining relation \rf{2} directly.	  
	\begin{proposition}\label{proposition 2} The center of $SO(2n+1)$ acts by multiplication on symmetric polynomials in $\gam_{2n,k}^2$. The center of $SO(2n)$ acts by multiplication of polynomials on $\gam^2_{2n-1,k}$, symmetric with respect to the permutations of the variables, and by powers of the monomial $\gam_{2n-1,1}\gam_{2n-,2}\cdots \gam_{2n-1,n}$.
	\end{proposition}
	This follows  from Harish-Chandra isomorphism, see e.g. \cite[Section 7.4]{D}. 
	
	Define the following automorphism of the space of meromorphic functions on $\gam_{kl}$ and $c$:
	\beq\label{7} \tau(\gam_{2k+1,j})=\gam_{2k+1,j},\qquad \tau(\gam_{2k,j})=-\gam_{2k,j},\qquad \tau(c)=-c\eeq
	\begin{lemma} We have the relations 
		\beq\label{8}\tau I_{k,k+1}=-I_{k,k+1}\tau \eeq
	\end{lemma}

	\setcounter{equation}{0}
	\section{Whittaker vectors}
	\subsection{Two chains of groups} Zhelobenko - Stern construction of the Gelfand - Tsetlin basis for orthogonal groups uses the chain of embeddings
	\beq\label{ch1}\ii_n:SO(N)\hookrightarrow SO(N+1)\eeq
	where the compact group $SO(N)$ is embedded into the compact group $SO(N+1)$ as the stabilizer of the vector $e_{N+1}$ so that the generators $I_{ij}$, $i,j\leq N$ of the Lie algebra $so(N)$ are identified with the corresponding generators of the Lie algebra $so(N+1)$.
	
	However, for the construction of Whittaker vectors in the related infinite dimensional representations of $so(N)$ in meromorphic functions we pass to another, noncompact  real form $SO(N,J)$ of the group $SO(N,\C)$ and use the chain of the corresponding Lie algebras compatible with the natural chain of Lie group $SO(N)$.
	Here
	$$J=\begin{pmatrix}0&0&...&0&1\\0&0&...&1&0\\&&...&&\\0&1&...&0&0\\1&0&...&0&0\end{pmatrix}$$
	The Lie algebra $so(N)$ is generated by elements $I_{ij}$, with the relation $I_{ji}=-I_{ij}$, so that the elements $I_{ij}$ with $i>j$ are chosen as a linear basis of Lie algebra $so(N)$. The Lie algebra $so(N,J)$, acting in the space with the basis $f_1,\ldots, f_N$ is generated by the elements 
	$$F_{ij}=f_{ij}-f_{\hat{j}\hat{i}},\qquad\text{where}\qquad f_{ij}(f_k)=\delta_{jk}f_i \qquad\text{and}\qquad\hat{i}=
	N+1-i$$ with the relation
	$F_{\hat{j},\hat{i}}=- F_{i,j}$
	so that the elements $F_{ij}$ with $i+j\leq N$ can be chosen as a linear basis of the Lie algebra $so(N,J)$.
	The elements $F_{ij}$ for $i<j$ form a positive nilpotent subalgebra, elements 
	$F_{kk}$ form a Cartan subalgebra.
	
	The chain of embedding
	\beq\label{ch2}\jj_n:SO(N,J)\hookrightarrow SO(N+1,J)\eeq
	is different.
	The group $SO(2n,J)$ is the stabilizer of the vector $f_{n+1}$ in the group $SO(2n+1,J)$, while the group $SO(2n-1,J)$ is the stabilizer of the element 
	$f_{n-1}+  (-1)^{n-1}  f_n$ in $SO(2n,J)$. Let us describe the maps 
	\beq\label{ch3} \ss_{N}:so(N)\to so(N,J) \eeq
	of {\bf complex} Lie algebras, which intertwine the embeddings \rf{ch1} and \rf{ch2}.
	On the level of bases of the vector space $\C^{2n}$ the map $\ss_{2n}$ corresponds to
	the transformation of  initial orthogonal basis 
	$e_1,\ldots, e_{2n}$
	of $\C^{2n}$ to the  defining basis $f_1,\ldots, f_{2n}$ of the form $J$,
	\beqq (f_1,\ldots, f_n,f_{n+1},\ldots, f_{2n}),\qquad (f_i,f_{\hat{j}})=\delta_{ij}\eeqq
by the relation 
\beq\label{fbasis}\begin{split} f_j=i^{j} \cdot \frac{ie_{2j-1}+e_{2j}}{\sqrt{2}},\qquad
	f_{\hat{j}}=i^{-j} \cdot \frac{-ie_{2j-1}+e_{2j}}{\sqrt{2}},\qquad j=1,,\ldots,n
	\end{split}\eeq 
	For the group $SO(2n+1)$ we transform the initial orthogonal basis 
	$e_1,\ldots, e_{2n+1}$
	of $\C^{2n+1}$ to the  defining basis $f_1,\ldots, f_{2n+1}$ of the form $J$, 
	\beqq (f_1,\ldots, f_n,f_{n+1},\ldots, f_{2n+1}),\qquad (f_i,f_{\hat{j}})=\delta_{ij}\eeqq 
	by the relation 
	\beqq
	\begin{split} f_j=i^{j} \cdot \frac{ie_{2j-1}+e_{2j}}{\sqrt{2}},\qquad
		f_{\hat{j}}=i^{-j} \cdot  \frac{-ie_{2j-1}+e_{2j}}{\sqrt{2}},\qquad j=1,,\ldots,n, \qquad f_{n+1}=e_{n+1}
	\end{split}\eeqq

	Correspondingly, the transformation formula from Lie algebra elements $I_{ij}$ to
	$F_{ij}$ are given by conjugation of the matrix $(I_{ij})$ be means of the corresponding transition matrix. In particular, we have the following expressions for the generators of Lie algebra $so(2n,J)$:
	\begin{align} \notag
		F_{j,j+1} = &\frac{1}{2}(I_{2j+1,2j} - I_{2j+2,2j-1}) +\frac{i}{2}(I_{2j+2,2j} + I_{2j+1,2j-1}) ,\\	    
		\notag  	F_{j+1,j} = &\frac{1}{2}(I_{2j+2,2j-1} - I_{2j+1,2j}) +\frac{i}{2}(I_{2j+2,2j} + I_{2j+1,2j-1},) \\	    
		\label{g20} 	F_{j,2n-j} =& (-1)^{j} \Big(  \frac{1}{2}(- I_{2j+1,2j} - I_{2j+2,2j-1}) +\frac{i}{2}(I_{2j+2,2j} - I_{2j+1,2j-1})  \Big), \\	   
		\notag	F_{2n-j,j} =& (-1)^{j} \Big( \frac{1}{2}(I_{2j+2,2j-1} + I_{2j+1,2j}) +\frac{i}{2}(I_{2j+2,2j} - I_{2j+1,2j-1})  \Big), && 
	\end{align}
	\begin{align} \label{g21}
		F_{j,j} = - i I_{2j,2j-1}
	\end{align}
	Here $j=1,\ldots, n-1$.
	Besides, instead of the use of the last simple root generator $F_{n-1,n}^{2n-1}$ of the Lie algebra $so(2n-1,J)$ it is convenient to use its image 
	\beqq \jj_{2n-1}(F_{n-1,n}^{2n-1})=\frac{1}{\sqrt{2}}\big(F_{n-1,n}^{2n}- (-1)^{n-1}  F_{n-1,n+1}^{2n}\big)
	\eeqq
	\beqq \jj_{2n-1}(F_{n,n-1}^{2n-1})=\frac{1}{\sqrt{2}}\big(F_{n,n-1}^{2n}-  (-1)^{n-1}  F_{n+1,n-1}^{2n}\big)
	\eeqq
	in the Lie algebra $so(2n)$. 
	
	Note that for the Lie algebra $so(n)$ the automorphism $\tau$ represents the longest element of the Weyl group,
	\beq\label{g20a} \tau F_{j,j+1}= F_{j+1,j}\tau. \eeq 
	\subsection{Right and left Whittaker vectors}
	For further convenience we denote by $\bgam_n$ the tuples of variables
	\beq\label{w0} \bgam_{2k}=\{\gam_{2k,1},\ldots \gam_{2k,k}\},\qquad \bgam_{2k-1}=\{\gam_{2k-1,1},\ldots \gam_{2k-1,k}\}\eeq
	and by $\hat{\bgam}_n$ the Gelfand - Tsetlin array
	\beq \label{w0a} \widehat{\bgam}_n= \begin{pmatrix}\bgam_n\\ \bgam_{n-1}\\ \vdots \\ \bgam_2\\ \bgam_1
	\end{pmatrix} \eeq 
	For any two sets
	\beqq \mathbf{X}= (x_1,\ldots,x_n),\qquad \mathbf{Y}=(y_1,\ldots, y_m)\eeqq
	set
	\beq \label{w1} s(\mathbf{X};\mathbf{Y})= \prod\limits_{x \in \mathbf{X},\, y \in \mathbf{Y}} (ic)^{\frac{x+y}{ic}} \cdot \Gamma \left(\frac{x+y}{ic} + \frac{1}{2}\right)\eeq
	With their help for any $k>0$ define the meromorphic functions $W_k^\pm$ and $V_k^\pm$ by the relation
	\beq \label{w2}\begin{split}
		W_k^\pm &= 
		e^{\mp\frac{\pi}{2c}\sum_j\gam_{2k,j}} 
		s(\pm\bgam_{2k-1},\bgam_{2k}) s(-\bgam_{2k},\pm\bgam_{2k+1}),\\
		V_k^\pm &= 
		e^{\mp\frac{\pi}{2c}\sum_j\gam_{2k,j}}
		s(\pm\bgam_{2k-1},\bgam_{2k}) .\\
	\end{split} 
	\eeq
	Set 
	\beq\label{w4}\begin{split} w_n&=
		 e^{-\frac{\pi n}{2c} \delta_{2n}}  W_1^+W_2^-\cdots W_{n-1}^{(-)^{n}}V_n^{(-)^{n+1}},\\
		w_n'&=
		e^{\frac{\pi}{c} \sum\limits_{k=1}^{n} \delta_{2k-1}} \tau(w_n) =  e^{-\frac{\pi n}{2c} \delta_{2n}}  W_1^-W_2^+\cdots W_{n-1}^{(-)^{n-1}}V_n^{(-)^{n}}
	\end{split}	\eeq
		Here 
	$$\delta_k=\sum_{j}\g_{kj}, \qquad \g_{kj}\in\bg_k. $$ 
	\begin{theorem}\label{theorem1}  The functions $w_n$ and $w'_n$ are left and right $SO(2n+1)$ Whittaker vectors:
		\begin{align}\label{t1} &F_{k,k+1}w_n=\frac{(-1)^{k+1}}{ic}w_n,\qquad k<n,&&F_{n,n+1}w_n=\frac{(-1)^{n  +1 }}{\sqrt{2}ic}w_n,\\
			\label{t2}	&F_{k+1,k}w'_n=\frac{(-1)^{k+1}}{ic}w'_n,\qquad k<n,
			&&F_{n+1,n}w'_n=\frac{(-1)^{n  +1 }}{\sqrt{2}ic}w'_n.
		\end{align}
	\end{theorem}
	{\bf Proof} is given in Appendix A
	
	Precise formulas for Whittaker vector and Whittaker function look better after the following change of variables
	\beq \label{g1}
	\g_{2k-1,j} = (-1)^{k+1} \nu_{2k-1,j},\qquad \g_{2k,j}=\gam_{2k,j}. 
	\eeq
	Then the Whittaker vectors can be written as
	\beq\label{g2}\begin{split}
		&w_n =  e^{-\frac{\pi n}{2c} \delta_{2n}}  e^{\sum\limits^{n}_{k=1} (-1)^{k} \delta_{2k}} \prod\limits^{n}_{k=1} s(\bg_{2k-1}, \bg_{2k}) \prod\limits^{n-1}_{k=1} s(-\bg_{2k}, -\bg_{2k+1}) \\ &
		w'_n =  e^{-\frac{\pi n}{2c} \delta_{2n}}  e^{\frac{\pi}{2c} \sum\limits^{n}_{k=1} (-1)^{k+1} \delta_{2k}} \prod\limits^{n}_{k=1} s(-\bg_{2k-1}, \bg_{2k}) \prod\limits^{n-1}_{k=1} s(-\bg_{2k},\bg_{2k+1} )
	\end{split}\eeq

	\subsection{Action of Cartan subalgebra}
	In Gelfand - Tsetlin representation of $gl_n$ see \cite{GKL}, the Cartan subalgebra acts by multiplication by linear functions on the variables $\g_{ij}$. It is not so for Gelfand-Tsetlin representations of $so(n)$, which we study here. However, the Cartan subalgebra of $so(n)$ acts in a similar way on Whittaker vectors.
	The Cartan subalgebra of $so(n)$ is generated by the elements $F_{kk}=-iI_{2k,2k-1}$. 
	\begin{proposition} 
		\begin{align}\label{c1} F_{kk}w_n &=\frac{(-1)^{k-1}}{ic} \left(\sum\limits_{j=1}^{k} \gam_{2k-1,j} + \sum\limits_{j=1}^{k-1} \gam_{2k-3,j} +(-1)^{k} (k-1)ic \right)w_n\\  F_{kk}w'_n & = \frac{(-1)^{k-1}}{ic} \left(\sum\limits_{j=1}^{k} \gam_{2k-1,j} + \sum\limits_{j=1}^{k-1} \gam_{2k-3,j} - (-1)^{k}i(k-1
			)c \right) w'_n\end{align}
	\end{proposition}
	In the variables $\g_{kj}$ the relations \rf{c1} look as follows 
	\begin{align}\label{c1a} F_{kk}w_n &=\frac{1}{ic} \left(\sum\limits_{j=1}^{k} \g_{2k-1,j} - \sum\limits_{j=1}^{k-1} \g_{2k-3,j} -(k-1) ic \right)w_n\\  F_{kk}w'_n & = \frac{1}{ic} \left(\sum\limits_{j=1}^{k} \g_{2k-1,j} - \sum\limits_{j=1}^{k-1} \g_{2k-3,j} + (k-1)ic \right) w'_n\end{align}
	
	 {\bf Proof} is given in Appendix B \cbk
	
	\setcounter{equation}{0}
	\section{Whittaker wave function}
	We construct Toda wave function  as a matrix coefficient \rf{g4} of elements of Cartan subgroup with respect to two Whittaker vectors. 
	\subsection{Invariant pairing}
	The invariant pairing looks the same in both $\gam$ and $\g$ variables.
	
	Define the functions $\tilde{\mu}(\bg_{2k})$ and  $\tilde{\mu}(\bg_{2k+1})$ by the relations
	\beq\label{m1} \begin{split}
		\tilde{\mu}(\bg_{2k})=e^{\frac{\pi}{c}\sum_j\g_{2k,j}}&\prod_r\Gamma^{-1}\Big(\frac{2\g_{2k,r}}{ic}\Big)\Gamma^{-1}\Big(\frac{-2\g_{2k,r}}{ic}\Big)\times\\ &\prod_{r\not=s}\Gamma^{-1}\Big(\frac{\g_{2k,r}-\g_{2k,s}}{ic}\Big)
		\prod_{r<s}\Gamma^{-1}\Big(\frac{\g_{2k,r}+\g_{2k,s}}{ic}\Big)\Gamma^{-1}\Big(\frac{-\g_{2k,r}-\g_{2k,s}}{ic}\Big),\\
		\tilde{\mu}(\bg_{2k+1})=&\prod_{r\not=s}\Gamma^{-1}\Big(\frac{\g_{2k+1,r}-\g_{2k+1,s}}{ic}\Big)
		\\	&\prod_{r<s}\Gamma^{-1}\Big(1+\!\frac{\g_{2k+1,r}+\g_{2k+1,s}}{ic}\Big)\Gamma^{-1}\Big(1-\!\frac{\g_{2k+1,r}+\g_{2k+1,s}}{ic}\Big),
	\end{split}\eeq
	that is 
	\beq\label{m1a} \begin{split}
		\tilde{\mu}(\bg_{2k})=e^{\frac{\pi}{c}\delta_{2k}}& \prod_{r<s}\left|\Gamma\Big(\frac{\g_{2k,r}-\g_{2k,s}}{ic}\Big)\right|^{-2}
		\left|\Gamma\Big(\frac{\g_{2k,r}+\g_{2k,s}}{ic}\big)\right|^{-2} \prod_r\left|\Gamma\Big(\frac{2\g_{2k,r}}{ic}\Big)\right|^{-2},\\
		\tilde{\mu}(\bg_{2k+1})=&\prod_{r<s}\left|\Gamma\Big(\frac{\g_{2k+1,r}-\g_{2k+1,s}}{ic}\Big)\right|^{-2}
		\left|\Gamma\Big(1+\!\frac{\g_{2k+1,r}+\g_{2k+1,s}}{ic}\Big)\right|^{-2}\,
	\end{split}\eeq
	and define a scalar product  on functions in $M$ as
	\beq\label{m2} (f,g)_{2n}=
	\int_{\R^{n^2}}\bar{f}(\bbg_{2n-1})g(\bbg_{2n-1})\tilde{\mu}(\bbg_{2n-1})d\bbg_{2n-1} \eeq
	where
	$$\tilde{\mu}(\bbg_n)=\prod_{k=1}^n \tilde{\mu}(\bg_k)$$
	Then 
	\begin{lemma}\label{lemma2}
		The operators $I_{2k+1,2k}$ and $	iI_{2k+2,2k+1}$ are skew symmetric with respect to the pairing \rf{m2}. In particular, operators $F_{kl}$ are 	 skew symmetric with respect to the pairing \rf{m2}.	
	\end{lemma} 
	\subsection{Integral formula}
	Due to \rf{g2} the product $\bar{w'}_n w_n$ looks as
	\beqq
	 \begin{split}\bar{w'}_n w_n= &  e^{-\frac{\pi n}{c} \delta_{2n}}  \prod\limits^{n-1}_{k=1} s(\bg_{2k-1}, \bg_{2k})  \bar{s}(-\bg_{2k-1}, \bg_{2k}) s(-\bg_{2k}, -\bg_{2k+1})\bar{s}(-\bg_{2k}, \bg_{2k+1})\times\\ &s(\bg_{2n-1},\bg_{2n})\bar{s}(-\bg_{2n-1},\bg_{2n})
	\end{split}\eeqq
	or, in terms of Gamma functions
	\beqq\begin{split}
		&w_{n} \cdot \overline{w_{n}'} =   
		 e^{-\frac{\pi}{c} \sum\limits^{n-1}_{k=1} \delta_{2k}}\cbk c^{\frac{1}{ic}\left( \sum\limits^{n-1}_{k=1} 2k(\delta_{2k-1} - \delta_{2k+1}) + 2n \delta_{2n-1} \right)} \\ &
		\times \left( \prod^{n-1}_{k=1} \Gamma \left(\frac{\pm \bg_{2k}+\bg_{2k-1}}{ic}+\frac{1}{2}\right) \Gamma \left(\frac{\pm \bg_{2k}-\bg_{2k+1}}{ic}+\frac{1}{2}\right) \right) \cdot \Gamma \left(\frac{\pm \bg_{2n}+\bg_{2n-1}}{ic}+\frac{1}{2}\right)
	\end{split}\eeqq
	Thus the Whittaker function
	\beq \label{g4}\Psi_{\bg_{2n}}=(w_n',\, e^{-\sum_{k=1}^nx_kF_{kk}}w_n)_{2n}\eeq
	is given by the integral
	\beq\label{g5}\begin{split}\Psi_{\bg_{2n}}(\bx_n)=&\int_{\R^{n^2}}
		\prod_{k=1}^{2n-1}{\mu}(\bg_k)d\bg_k\,
		\cdot c^{\frac{2}{ic}\sum_{k=1}^n\delta_{2k-1}}
		e^{\frac{1}{ic}\sum_{k=1}^n(\delta_{2k-3}-\delta_{2k-1}+(k-1)ic)x_k}\times\\
		&\prod_{k=1}^{n-1}\prod_{i,j=1}^k\prod_{l=1}^{k+1}\Gamma\left(\frac{\pm\g_{2k,i}+\g_{2k-1,j}}{ic}+\frac{1}{2} \right)\Gamma\left(\frac{\pm\g_{2k,i}-\g_{2k+1,l}}{ic}+\frac{1}{2} \right)\times\\ &\prod_{i,j=1}^n\Gamma\left(\frac{\pm\g_{2n,i}+\g_{2n-1,j}}{ic}+\frac{1}{2} \right)
	\end{split}	\eeq
	Here
	\beq\label{m1b} \begin{split}
		{\mu}(\bg_{2k})=e^{-\frac{\pi}{c}\delta_{2k}}\tilde{\mu}(\bg_{2k})&= \prod_{r<s}\left|\Gamma\Big(\frac{\g_{2k,r}-\g_{2k,s}}{ic}\Big)\right|^{-2}
		\left|\Gamma\Big(\frac{\g_{2k,r}+\g_{2k,s}}{ic}\big)\right|^{-2} \prod_r\left|\Gamma\Big(\frac{2\g_{2k,r}}{ic}\Big)\right|^{-2},\\
		{\mu}(\bg_{2k+1})=\tilde{\mu}(\bg_{2k+1})&=\prod_{r<s}\left|\Gamma\Big(\frac{\g_{2k+1,r}-\g_{2k+1,s}}{ic}\Big)\right|^{-2}
		\left|\Gamma\Big(1+\!\frac{\g_{2k+1,r}+\g_{2k+1,s}}{ic}\Big)\right|^{-2}\,
	\end{split}\eeq
	The measure functions ${\mu}(\bg_j)$ 
	do not contain exponential factors.
	
	The convergence of the integral \rf{g5} can be proved by the arguments given in \cite[Appendix A]{IS}.
	 Namely, let us complete the sequence 
	 \beqq \bg_1=\{\g_{11}\}, \bg_2=\{\g_{21}\},\bg_3=\{\g_{31},\g_{32}\},\ldots \bg_{2n}=\{\g_{2n,1},\ldots \g_{2n,n}\},\eeqq
	 to the sequence
	 \beqq \bg'_1=\{\g'_{1,-1},\g'_{11}\}, \bg'_2=\{\g'_{2,-1},\g'_{2,0},\g'_{21}\},,\ldots \bg'_{2n}=\{\g_{2n,-n},\ldots \g_{2n,n}\},\eeqq
	where
	\beqq\g'_{n,k}=\g_{n,k}, \ \g'_{n,-k}=-\g_{n,k}\qquad\text{for}\qquad k>0,\qquad\text{and}\qquad \g'_{2m,0}=0\eeqq
	as   is customary in the representation theory of orthogonal groups. Then the inequality \cite[34]{IS} applied to this sequence is transformed to the bound
	\beq\label{g14a}\begin{split}&\sum_{k=1}^{2n}\sum_{i,j}|\pm \g_{2k,i}-\g_{2k\pm 1,j}|-\sum_{r=1}^{2n-1}\sum_{i\not=j}|\pm\g_{r,i}-\g_{r,j}|-2\sum_{k=1}^{n-1}\sum_i|\g_{2k,i}|\geq \\
	&	C(\bg_{2n})+\frac{2}{n}\sum_{r=1}^{2n-1}\sum_i|\g_{r,i}|,
				\end{split}\eeq
			where the constant $C(\bg_{2n})$ depends on the values of $\g_{2n,i}$.
			Due to the asymptotics of the Gamma function $$\Gamma(ix)\sim  e^{-\pi|x|/{2}}$$ in imaginary direction we observe that the integrand of \rf{g5} can be bounded by \rf{g14a} as
			\beq  C'(\bg_{2n}) e^{-\frac{\pi-\ve}{cn}\sum_{r=1}^{2n-1}\sum_i|\g_{r,i}|} \eeq
			for any small positive $\ve>0$ and a proper positive constant  $C'(\bg_{2n})$, which implies absolute convergence of the integral \rf{g5}
	\subsection{Toda equation}
	Denote by $H_n^B$ the Toda hamiltonian 
	\beq\label{g16} H_{B_n}= \sum_{k=1}^n\left(-\frac{\d^2}{\d x_k^2}+(2n-2k+1)\frac{\d}{\d x_k}\right) +\sum_{k=1}^{n-1}\frac{2}{c^2}e^{x_k-x_{k+1}}+\frac{1}{c^2}e^{x_n}\eeq
	\begin{theorem}\label{theorem2} The Whittaker function \rf{g5} is a wave function for
		$B_n$ Toda Hamiltonian:
		\beq\label{g17}
		H_{B_n}\Psi_{\bg_{2n}}(\bx_n)=\frac{1}{c^2}\sum_{i=1}^n\g_{2n,i}^2+\frac{n(2n-1)(2n+1)}{12}\eeq
	\end{theorem}
	{\bf Proof} is a standard game with the matrix element
	\beqq (w_n'\, L_{2n+1} \,e^{-\sum_k x_kF_{kk}} w_n)_{2n} \eeqq
	where 
	\beqq 
	L_{2n+1}=\sum_{i,j=1}^{2n+1}I_{ij}^2=\sum_{i,j=1}^{2n+1}F_{ij}F_{ji}\eeqq
	is Laplace operator of $SO(2n+1)$.
	For this one should also know the eigenvalue of $L$ in the representation in meromorphic functions. But it is known from the  theory of highest weight representations of $so(2n+1)$. It gives us the eigenvalue
	\beqq  \sum_{i=1}^{n}m_{2n,i}^2- (\rho,\rho)\eeqq
	where 
	\beqq \rho=(n-\frac{1}{2},\ldots, \frac{1}{2})\eeqq
	is a half sum of positive roots for $so(2n+1)$.
	Thus the eigenvalue is equal to
	\beqq -\frac{1}{c^2} \sum_{i=1}^{n}\g_{2n,i}^2-\frac{n(2n-1)(2n+1)}{12}\eeqq
	\hfill{$\Box$}
	
{\bf Remark}.	Note that the function
	\beqq \tilde{\Psi}_{\bg_{2n}}(\bx_n)=e^{-(\rho,\bx_n)}{\Psi}_{\bg_{2n}}(\bx_n)=
	e^{-\sum\limits^{n}_{k=1} (n-k+\frac{1}{2})x_{k}}{\Psi}_{\bg_{2n}}(\bx_n)\eeqq
	is the solution of more familiar spectral problem
	\beqq 
	 \left(-\sum_{k=1}^n\frac{\d^2}{\d x_k^2} +\sum_{k=1}^{n-1}\frac{2}{c^2}e^{x_k-x_{k+1}}+\frac{1}{c^2}e^{x_n}\right)\tilde{\Psi}_{\bg_{2n}}(\bx_n)=
	\frac{1}{c^2}\left(\sum_{i=1}^n\g_{2n,i}^2\right)\tilde{\Psi}_{\bg_{2n}}(\bx_n) \eeqq

	\subsection{Iterative procedures}

	{\bf 1}. The integral \rf{g5} can be formulated as an iterative integral presentation of the Whittaker wave function,
	\beq\label{g6}\begin{split}\Psi_{\bg_{2n}}(x_1,\ldots,x_n)=\int_{\R^{2n-1}}{\mu}(\bg_{2n-1})d\bg_{2n-1}\,
		{\mu}(\bg_{2n-2})d\bg_{2n-2}\\
		\prod_{i=1}^{n-1}\prod_{l=1}^{n}\Gamma\left(\frac{\pm\g_{2n-2,i}-\g_{2n-1,l}}{ic}+\frac{1}{2} \right) \prod_{i,j=1}^n\Gamma\left(\frac{\pm\g_{2n,i}+\g_{2n-1,j}}{ic}+\frac{1}{2} \right)\times\\c^{\frac{2}{ic}\delta_{2n-1}}e^{-\frac{\delta_{2n-1}}{ic}x_n+((n-1)  +  \sum_{k=1}^{n-2}k)x_n}\Psi_{\bg_{2n-2}}(x_1-x_n,\ldots, x_{n-1}-x_n)
	\end{split}\eeq
	or
	\beq\label{g7}\Psi_{\bg_{2n}}(x_1,\ldots,x_n)=\Lambda_n(x_n)\left(\Psi_{\bg_{2n-2}}(x_1-x_n,\ldots,x_{n-1}-x_n)\right)\eeq
	where $\Lambda(x_n)$ is an integral operator
	\beq \label{g8} (\Lambda_n(x)f)(\bg_{2n})=\int_{\R^{n-1}}K(\bg_{2n};\bg_{2n-2}|x)f(\bg_{2n-2})\mu(\bg_{2n-2})d\bg_{2n-2} \eeq
	with the kernel
	\beq\label{g9}\begin{split} K(\bg_{2n};\bg_{2n-2}|x)=e^{\frac{n(n-1)}{2}x}\int_{\R^{n}} {\mu}(\bg_{2n-1})d\bg_{2n-1}c^{\frac{2}{ic}\delta_{2n-1}}e^{-\frac{\delta_{2n-1}}{ic}x}\\	\prod_{i=1}^{n-1}\prod_{l=1}^{n}\Gamma\left(\frac{\pm\g_{2n-2,i}-\g_{2n-1,l}}{ic}+\frac{1}{2} \right) \prod_{i,j=1}^n\Gamma\left(\frac{\pm\g_{2n,i}+\g_{2n-1,j}}{ic}+\frac{1}{2} \right)
	\end{split}\eeq 
	
	{\bf 2}. Another recurrent procedure uses the observation that the above construction of the Whittaker vectors and Whittaker functions, restricted to $SO(2n)$ produced actually $GL_n$ Whittaker vectors and functions. It can be seen from the relation \rf{A2}, \rf{A2a}. In this way we arrive by using Gustafson integrals \cite{Gust} to Iorgov - Shadura formula \cite{IS}, which expresses the $B_n$  Toda wave function via $A_n$ Toda wave function and contains in total twice less integrals.
	
	The restriction $\Psi_{\bg_{2n-1}}(\bx_n)$ of the wave function \rf{g5} to $SO(2n)$ is given by the integral 
	\beq\label{g10}\begin{split}\Psi_{\bg_{2n-1}}(\bx_n)=&\int_{\R^{n^2-n}}
		\prod_{k=1}^{2n-2}{\mu}(\bg_k)d\bg_k\,
		\cdot c^{\frac{2}{ic}\sum_{k=1}^{n  - 1 } \delta_{2k-1}}
		e^{\frac{1}{ic}\sum_{k=1}^n(\delta_{2k-3}-\delta_{2k-1}+(k-1)ic)x_k}\times\\
		&\prod_{k=1}^{n-1}\prod_{i,j=1}^k\prod_{l=1}^{k+1}\Gamma\left(\frac{\pm\g_{2k,i}+\g_{2k-1,j}}{ic}+\frac{1}{2} \right)\Gamma\left(\frac{\pm\g_{2k,i}-\g_{2k+1,l}}{ic}+\frac{1}{2} \right)
	\end{split}	\eeq
	and the functions $\Psi_{\bg_{2n}}(\bx_n)$ and $\Psi_{\bg_{2n-1}}(\bx_n)$ are related as follows
	\beq\label{g11} \Psi_{\bg_{2n}}(\bx_n)=\int_{\R^n}
	{\mu}(\bg_{2n-1})d\bg_{2n-1}c^{\frac{2}{ic}\delta_{2n-1}}
	\prod_{i,j=1}^n\Gamma\left(\frac{\pm\g_{2n,i}+\g_{2n-1,j}}{ic}+\frac{1}{2} \right) \Psi_{\bg_{2n-1}}(\bx_n)
	\eeq
	 For each $k$ the integral over $\bg_{2k}$ in \rf{g10} can be explicitly calculated by means of the degenerate $B_n$ Gustafson integral
	\beq\label{gu1}\begin{split}
		\frac{1}{(2\pi)^n}\int\limits_{\R^n}
		\frac{\prod\limits_{i=1}^{2n+1}\prod\limits_{j=1}^n\Gamma(a_i+iz_j))
			\Gamma(a_i-iz_j))}{\prod\limits_{1\leq i<j\leq n}\big|\Gamma(i(z_i+z_j))
			\Gamma(i(z_i-z_j))\big|^2
			\prod\limits_{j=1}^n\big|\Gamma(2iz_j)\big|^2}
		d\bx_n =
		{n!2^n\prod\limits_{1 \leq i < j \leq 2n+1}\Gamma(a_i+a_j)}
		\;,
	\end{split}\eeq
	where all $a_i$ are assumed to have a positive real part. The integral \rf{gu1}  is a limiting case  
	$$a_{2n+2}=\ve+iL,\qquad  L\to\infty, \qquad \ve\to +0$$
	of a general $B_n$ Gustafson integral \cite{Gust}
	\beqq
	\frac{1}{(2\pi)^n}\int\limits_{\R^n}
	\frac{\prod\limits_{i=1}^{2n+2}\prod\limits_{j=1}^n\Gamma(a_i+iz_j))
		\Gamma(a_i-iz_j))}{\prod\limits_{1\leq i<j\leq n}\big|\Gamma(i(z_i+z_j))
		\Gamma(i(z_i-z_j))\big|^2
		\prod\limits_{j=1}^n\big|\Gamma(2iz_j)\big|^2}
	d\bx_n =
	\frac{n!2^n\prod\limits_{1\leq i<j\leq 2n+2}\Gamma(a_i+a_j)}
	{\Gamma\Big(\sum\limits_{i=1}^{2n+2}a_i\Big)}.
	\eeqq
	Using 
	\beqq a_1=\frac{\g_{2k-1,1}}{ic}+\frac{1}{2},\ldots, a_k=\frac{\g_{2k-1,k}}{ic}+\frac{1}{2},
	\qquad a_{k+1}= -\frac{\g_{2k+1,1}}{ic}+\frac{1}{2},\ldots a_{2k+1}=-\frac{\g_{2k+1,k+1}}{ic}+\frac{1}{2}, \eeqq
	we get 
	\beq\label{g12}\begin{split}
		& \int\limits_{\mathbb{R}^{k}} \frac{\prod\limits^{k}_{i=1} \left(\prod\limits^{k}_{j=1} \Gamma \left(\frac{\pm \g_{2k,i}+\g_{2k-1,j}}{ic}+\frac{1}{2}\right) \prod\limits^{k+1}_{j=1} \Gamma \left(\frac{\pm \g_{2k,i}-\g_{2k+1,j}}{ic}+\frac{1}{2}\right) \right)}{\prod\limits_{r<s}\left|\Gamma\Big(\frac{\g_{2k,r}-\g_{2k,s}}{ic}\Big)\right|^{2}	\left|\Gamma\Big(\frac{\g_{2k,r}+\g_{2k,s}}{ic}\big)\right|^{2} \prod\limits_{r}\left|\Gamma\Big(\frac{2\g_{2k,r}}{ic}\Big)\right|^{2}} d\g_{2k,1} \ldots d\g_{2k,k} = \\ &
		= c^{k} (2 \pi)^{k} 2^{k} \cdot k! \cdot \prod\limits_{r<s} \Gamma \left(\frac{ \g_{2k-1,r}+\g_{2k-1,s}}{ic}+1 \right) \cdot \prod\limits_{i,j} \Gamma \left(\frac{ \g_{2k-1,i}-\g_{2k+1,j}}{ic}+1 \right) \\ &
		\qquad \qquad \times \prod\limits_{r<s} \Gamma \left(-\frac{ \g_{2k+1,r}+\g_{2k+1,s}}{ic}+1 \right)
	\end{split}\eeq
	and 
	\beq \label{g13}
	\Psi_{\bg_{2n-1}}(\bx_n)=d_{n} \cdot c^{\frac{n+1}{ic} \delta_{2n-1}}\prod_{1\leq r <s\leq n}\Gamma \left(1-\frac{ \g_{2n-1,r}+\g_{2n-1,s}}{ic} \right)\Phi_{\bg_{2n-1}}(\bx_n)\eeq
	where 
	\beq\label{g14}\begin{split}
		\Phi_{\bg_{2n-1}} (\bx_n )= \int\limits_{\mathbb{R}^{n(n-1)/2}}
		e^{\frac{1}{ic} \sum\limits^{n}_{k=1} (\delta_{2k-3} - \delta_{2k-1} + (k-1)ic)x_{k}} d\g_{1} \ldots d\g_{2n-3}\times \\
		\frac{\prod\limits^{n-1}_{k=1} \prod\limits^{k}_{i=1} \prod\limits^{k+1}_{j=1} c^{\frac{ \g_{2k-1,i}-\g_{2k+1,j}}{ic}} \Gamma \left(\frac{ \g_{2k-1,i}-\g_{2k+1,j}}{ic}+1 \right)}{\prod\limits^{n-1}_{k=1} \prod\limits_{r<s} \left|\Gamma \left(\frac{ \g_{2 k -1,r}-\g_{2 k -1,s}}{ic}\right)\right|^{2}} 
	\end{split}\eeq
	and 
	\beq\label{g15} d_n=\prod_{k=1}^{n-1} c^{k} (2 \pi)^{k} 2^{k} \cdot k!\eeq
	Both functions 	$\Psi_{\bx_n} (\g_{2n-1} )$ and $\Phi_{\bg_{2n-1}}(\bx_n )$
	are solutions of $GL(n)$ Toda equations. In order to compare the final results with  Iorgov - Shadura formula \cite{IS} we perform the change of integration variables
	$$ \g_{2k-1,j}\to \g_{2k-1,j}-\left(n-k+\frac{1}{2}\right)ic.$$
	Then 
	\beq\label{g16a}\begin{split}
		\Phi_{\bg_{2n-1}} (\bx_n )= e^{\sum\limits^{n}_{k=1} (n-k+\frac{1}{2})x_{k}}  &\int\limits_{C_n} e^{\frac{1}{ic} \sum\limits^{n}_{k=1} (\delta_{2k-3} - \delta_{2k-1})x_{k}} d\g_{1} \ldots d\g_{2n-3}\times \\	 	
		&\frac{\prod\limits^{n-1}_{k=1} \prod\limits^{k}_{i=1} \prod\limits^{k+1}_{j=1} c^{\frac{ \g_{2k-1,i}-\g_{2k+1,j}}{ic}} \Gamma \left(\frac{ \g_{2k-1,i}-\g_{2k+1,j}}{ic} \right)}{\prod\limits^{n-1}_{k=1} \prod\limits_{r<s} \left|\Gamma \left(\frac{ \g_{2 k -1,r}-\g_{2 k -1,s}}{ic}\right)\right|^{2}} 
	\end{split}\eeq
	and
	\beq\label{g17a}\begin{split}
		 \Psi_{\bg_{2n}} (\bx_n )=d_n& \int\limits_{C'_n} \frac{\prod\limits^{n}_{i=1}\prod\limits^{n}_{j=1} \Gamma \left(\frac{ \g_{2n-1,i}+\g_{2n,j}}{ic}\right) \Gamma \left(\frac{ \g_{2n-1,i}-\g_{2n,j}}{ic} \right)}{\prod\limits_{r<s} \Gamma \left(\frac{ \g_{2n-1,r}+\g_{2n-1,s}}{ic} \right)  \prod\limits_{r<s} \left|\Gamma \left(\frac{ \g_{2n-1,r}-\g_{2n-1,s}}{ic}\right)\right|^{2} }\times \\ & c^{\frac{n+1}{ic} \delta_{2n-1}} \Phi_{\bg_{2n-1}} (\bx_n) d\g_{2n-1}
	\end{split}\eeq
Here the contour $C_n$ in \rf{g16a} is a deformation of $\mathbb{R}^{n(n-1)/2}$, such that the integration over the variable $\g_{2k-1}$ is performed in such a way that the singularity of $\Gamma \left(\frac{ \g_{2k-1,i}-\g_{2k+1,j}}{ic} \right)$ is under the line of integration and the  singularity of $\Gamma \left(\frac{ \g_{2k-3,i}-\g_{2k-1,j}}{ic} \right)$ is above the line of integration  over $\g_{2k-1}$. The contour $C'_n$ in \rf{g17a} is a deformation of $\mathbb{R}^{n}$ where the singularities of the nominators are under the contours of integrations over all variables.
 
	The relations \rf{g16a} - \rf{g17a} are in accordance  with Iorgov - Shadura formula \cite[(26),(27)]{IS}. More precisely, in Iorgov - Shadura description the boundary wall corresponds to the first coordinate $x_1$, while we work with the boundary wall related to the last coordinate. One can observe the coincidence of  formulas after the change of variables
	\beqq x_k\to -x_{n+1-k} \eeqq
	 in their (or in ours) formulas.

	\setcounter{equation}{0}
		\section{Examples} 
		{$\ \mathbf{ n=1}$}. For one particle the system and the wave function coincide with that of $sl(2)$ Toda
		system:
		\beqq \Psi_{\g_{21}}(x_1)=\int_\R d\g_{11} c^{\frac{2}{ic}\g_{11}}e^{-\frac{1}{ic}\g_{11} x_1}
		\Gamma\left(\frac{\pm\g_{21}+\g_{11}}{ic}+\frac{1}{2}\right)\eeqq
	{$\ \mathbf{ n=2}$}. The wave function for $B_2$ Toda system is given by $4$  fold integral, which can be reduced to three fold by using Gustafson integral,
	\beqq\begin{split}\Psi_{\g_{41},\g_{42}}(x_1,x_2)=\int_{\R^4}& d\g_{11}
	\frac{d \g_{21}}{\Big|\Gamma\Big(\frac{2\g_{21}}{ic} \Big)\Big|^2} 
	\frac{d \g_{31} d\g_{32}}{\Big|\Gamma\Big(\frac{\g_{31}-\g_{32}}{ic} \Big)
		\Gamma\Big(1+\frac{\g_{31}+\g_{32}}{ic} \Big)\Big|^2} \times\\  & c^{\frac{2}{ic}(\g_{11}+\g_{31}+\g_{32})}
	e^{\frac{1}{ic}\left( -\g_{11}x_1+(\g_{11}-\g_{31}-\g_{32}+ic)x_2\right)}\times\\	\Gamma\left(\frac{\pm\g_{21}+\g_{11}}{ic}+\frac{1}{2}\right)&
	\Gamma\left(\frac{\pm\g_{21}-\g_{31}}{ic}+\frac{1}{2}\right)
	\Gamma\left(\frac{\pm\g_{21}-\g_{32}}{ic}+\frac{1}{2}\right)
	\prod_{i,j=1}^2\Gamma\left(\frac{\pm\g_{4i}+\g_{3j}}{ic}+\frac{1}{2}\right)
	\end{split}\eeqq
that is
  \beqq\begin{split}&\Psi_{\g_{41},\g_{42}}(x_1,x_2)=\\ &\int_{\R^3} 
  	\frac{d \g_{21}}{\Big|\Gamma\Big(\frac{2\g_{21}}{ic} \Big)\Big|^2} 
  	\frac{d \g_{31} d\g_{32}}{\Big|\Gamma\Big(\frac{\g_{31}-\g_{32}}{ic} \Big)
  		\Gamma\Big(1+\frac{\g_{31}+\g_{32}}{ic} \Big)\Big|^2} c^{\frac{2}{ic}(\g_{31}+\g_{32})}
  	e^{\frac{1}{ic}\left( (-\g_{31}-\g_{32}+ic)x_2\right)}\\	&
  	\Gamma\left(\frac{\pm\g_{21}-\g_{31}}{ic}+\frac{1}{2}\right)
  	\Gamma\left(\frac{\pm\g_{21}-\g_{32}}{ic}+\frac{1}{2}\right)
  	\prod_{i,j=1}^2\Gamma\left(\frac{\pm\g_{4i}+\g_{3j}}{ic}+\frac{1}{2}\right)
  	 \Psi_{\g_{21}}(x_1-x_2)
  \end{split}\eeqq
	or
	\beqq\begin{split}\Psi_{\g_{41},\g_{42}}(x_1,x_2)=\int_{\R^2+\ve}
\frac{\prod_{i,j=1}^2\Gamma\left(\frac{\pm\g_{4i}+\g_{3j}}{ic}\right)}{ \Gamma\left(\frac{\g_{31}+\g_{32}}{ic}\right) \left|\Gamma \left(\frac{ \g_{31}-\g_{32}}{ic}\right)\right|^{2} }c^{\frac{3}{ic}(\g_{31}+\g_{32})}\Phi_{\g_{31},\g_{32}}(x_1,x_2) d\g_{31}d\g_{32} 		\end{split}\eeqq
	where $\Phi_{\g_{31},\g_{32}}(x_1,x_2)$ is the wave function of $GL(2)$ Toda system 
	\beqq\begin{split}
		\Phi_{\g_{31},\g_{32}}(x_1,x_2)= e^{\frac{3x_1+x_2}{2}}&\int_{\R+\ve} d\g_{11} e^{\frac{1}{ic}(-\g_{11}x_1+(\g_{11}-\g_{31}-\g_{32})x_2)}c^{\frac{1}{ic}(2\g_{11}-\g_{31}-\g_{32})}\\& \Gamma\left(\frac{\g_{11}-\g_{31}}{ic}\right) \Gamma\left(\frac{\g_{11}-\g_{32}}{ic}\right)
		\end{split}\eeqq
	\setcounter{equation}{0}
	\appendix
	\renewcommand*{\thesection}{A}
\section*{Acknowlegements} The authors thank N.Belousov, S.Derkachov and S.Kharchev for their interest in the work and for fruitul discussions. 
The work of S. Kh. is an output of a research project implemented as part of the Basic Research Program at the National Research University Higher School of Economics (HSE University).
\setcounter{equation}{0}
\section{Calculations of Whittaker vector}
In this section we prove Theorem \ref{theorem1}.
First we note that the relation \rf{t2} follows from \rf{t1} by using the automorphism $\tau$.  
 Proof of the  equality \rf{t1} reduces to  check of the following equalities, where  $k=1,\ldots n-1$ in the relations \rf{A2a} and \rf{A2b}    
\begin{align} \label{A2} 
	& (F_{k,k+1} -  (-1)^{k}  F_{k,2n-k})w_{n} = (I_{2k+1,2k} + iI_{2k+1,2k-1})w_{n}= \frac{ (-1)^{k+1} }{ic}w_{n}  \\ \label{A2a} &
	(F_{k,k+1} +   (-1)^{k} \cbk F_{k,2n-k})w_{n} = (- I_{2k+2,2k-1} + i I_{2k+2,2k})w_{n}=\frac{ (-1)^{k+1} }{ic}w_{n}  \\ \label{A2b} &
	F_{n,n+1}w_{n} = \frac{1}{\sqrt{2}} (I_{2n+1,2n} + iI_{2n+1,2n-1})w_{n}= \frac{ (-1)^{n+1} }{ic\sqrt{2}}w_{n}
\end{align}
Besides, the relation \rf{A2b} is a particular case of \rf{A2}. So we have to prove the relations  \rf{A2} and \rf{A2a}.
\medskip

The proof of \rf{A2}--\rf{A2a} requires certain calculations. For their visualization we introduce some intermediate notations. First rewrite the operators \rf{5} and \rf{6} as
\beq\label{A3} \begin{split}
	I_{2k+1,2k}=\sum \limits_{\ve = \pm 1} \sum\limits_{j=1}^{k} P^{\ve}_{kj},\qquad
	I_{2k+2,2k+1}= \sum\limits_{j=1}^{k} Q_{kj} + \sum\limits_{j=1}^{k} R_{kj} + T_k
\end{split}
\eeq 
where 
\beqq\begin{split}
	P^{\ve}_{k,j} =& - \frac{1}{ic}  \frac{\prod\limits_{r=1}^{k-1} (\gam_{2k-1,j} + \ve (\gam_{2k-2,r} + \frac{ic}{2})) \prod\limits_{r=1}^{k} (\gam_{2k-1,j} - \ve (\gam_{2k,r} - \frac{ic}{2}))}{2\prod\limits_{r \neq j} (\gam_{2k-1,j} - \gam_{2k-1,r})(\gam_{2k-1,j} + \gam_{2k-1,r} + \ve ic)} e^{\ve ic \partial_{\gam_{2k-1,j}}},\\
	Q_{k,j} = & \frac{1}{ic}  \frac{\prod\limits_{r=1}^{k+1} ((\gam_{2k,j}+\frac{ic}{2})^2 - \gam_{2k+1,r}^2)}{2\gam_{2k,j} (\gam_{2k,j}+\frac{ic}{2}) \prod\limits_{r \neq j} (\gam_{2k,j}^2 - \gam_{2k,r}^2)} e^{ ic \partial_{\gam_{2k,j}}} \\
	R_{k,j} =& \frac{1}{ic} \frac{\prod\limits_{r=1}^{k} ((\gam_{2k,j}-\frac{ic}{2})^2 - \gam_{2k-1,r}^2)}{2\gam_{2k,j} (\gam_{2k,j}-\frac{ic}{2}) \prod\limits_{r \neq j} (\gam_{2k,j}^2 - \gam_{2k,r}^2)} e^{-ic \partial_{\gam_{2k,j}}} \\
	T_{k} =& \frac{1}{c} \frac{\prod\limits_{r=1}^{k} \gam_{2k-1,r} \prod\limits_{r=1}^{k+1} \gam_{2k+1,r}}{\prod\limits_{r=1}^{k} (\gam_{2k,r}^2 + \frac{c^2}{4})}
\end{split} \eeqq
Set also for $\ve,\delta=\pm 1$
\beq\label{A4}
J^{\ve}_{k,\delta,j} =  \sum\limits_{s=1}^{k} \frac{\ve ic}{\gam_{2k+\delta,j} - \ve (\gam_{2k,s}+\frac{ic}{2})} Q_{k,s} + \sum\limits_{s=1}^{k} \frac{\ve ic}{\gam_{2k + \delta,j} + \ve (\gam_{2k,s}-\frac{ic}{2})} R_{k,s} + \frac{\ve ic}{\gam_{2k + \delta,j}} T_{k}. 
\eeq
In these notation we can present the following expressions for other generators of Lie algebra $so(n)$ needed in the equations on Whittaker vectors. They can be checked by straightforward calculations. 
\begin{lemma}\label{lemmaA1} We  have the following relations      
	\beq\label{A5}\begin{split}
		&I_{2k+1,2k-1} = \sum \limits_{\ve = \pm 1} \sum\limits_{j=1}^{k} P^{\ve}_{kj} J_{k-1,1,j}^{\ve} \qquad k = 1,\ldots,n \\ &
		I_{2k+2,2k} = - \sum \limits_{\ve = \pm 1} \sum\limits_{j=1}^{k} P^{\ve}_{kj} J_{k,-1,j}^{\ve} \qquad k = 1,\ldots,n-1 \\ &
		I_{2k+2,2k-1} = - \sum \limits_{\ve = \pm 1} \sum\limits_{j=1}^{k} P^{\ve}_{kj} J_{k,-1,j}^{\ve} J_{k-1,1,j}^{\ve} \qquad k = 1,\ldots,n-1 \\ &
	\end{split} \eeq
\end{lemma}
For more brevity in futher formulas we denote
$$\theta_k=(-1)^k,\qquad \theta_{k+1}=(-1)^{k+1}.$$
The following statement is one of the important steps in the proof of Theorem 1.
\begin{lemma}\label{lemmaA2} For any 
	$k=0,\ldots,n-1$ ($\delta = 1$ when $k=0$) we have the relations
	\begin{align} \label{A6}
		&J^{\tt_k}_{k,\delta,j} w_{n} = i w_{n} \\ \label{A7}&
		J^{\tt_{k+1}}_{k,\delta,j} w_{n} = i \left( 1 - \frac{\prod\limits_{r=1}^{k+1} (\gam_{2k+\delta,j} + \gam_{2k+1,r}) \prod\limits_{r=1}^{k} (\gam_{2k+\delta,j} + \gam_{2k-1,r})}{\gam_{2k+\delta,j} \prod\limits_{r=1}^{k} ((\gam_{2k+\delta,j}+(-1)^{k}\frac{ic}{2})^2 - \gam_{2k,r}^2)} \right) w_{n} \qquad \delta = \pm 1
	\end{align}
\end{lemma}
{\bf Proof .} Note that the operator $J^{\tt_{k}}_{k,\delta ,j}$ containes only shifts of variables  $\gam_{2k,s}$, and thus only the factor  $W_{k}^{\tt_{k+1}}$ can change. Let us check the equality \rf{A6}.
We then have 
\beqq \begin{split} 
	&J^{\tt_{k}}_{k,\delta,j} w_{n} = \\ \bigg(\sum\limits_{s=1}^{k}& \frac{\tt_{k} ic}{\gam_{2k+\delta,j} - \tt_k  (\gam_{2k,s}+\frac{ic}{2})} Q_{k,s} + \sum\limits_{s=1}^{k} \frac{\tt_{k} ic}{\gam_{2k+\delta,j} + \tt_k  (\gam_{2k,s}-\frac{ic}{2})} R_{k,s} + \frac{\tt_k ic}{\gam_{2k + \delta,j}} T_{k} \bigg) w_{n} \\
	=&  
	\Bigg( \sum\limits_{s=1}^{k} \frac{\tt_k ic}{\gam_{2k+\delta,j} - \tt_k (\gam_{2k,s}+\frac{ic}{2})} \cdot \frac{1}{ic} \cdot \frac{\prod\limits_{r=1}^{k+1} ((\gam_{2k,s}+\frac{ic}{2})^2 - \gam_{2k+1,r}^2)}{2\gam_{2k,s} (\gam_{2k,s}+\frac{ic}{2}) \prod\limits_{r \neq s} (\gam_{2k,s}^2 - \gam_{2k,r}^2)} \times \\  
	i^{\tt_k} \times & \frac{\prod\limits_{r=1}^{k} (\gam_{2k,s} - \tt_k \gam_{2k-1,r} + \frac{ic}{2})}{\prod\limits_{r=1}^{k+1} (-\gam_{2k,s} -\tt_k\gam_{2k+1,r} - \frac{ic}{2})} +
	\sum\limits_{s=1}^{k} \frac{\tt_k ic}{\gam_{2k+\delta,j} + \tt_k (\gam_{2k,s}-\frac{ic}{2})} \cdot \frac{1}{ic} \times\\ &\frac{\prod\limits_{r=1}^{k+1} ((\gam_{2k,s}-\frac{ic}{2})^2 - \gam_{2k+1,r}^2)}{2\gam_{2k,s} (\gam_{2k,s}-\frac{ic}{2}) \prod\limits_{r \neq s} (\gam_{2k,s}^2 - \gam_{2k,r}^2)} \cdot
	i^{-\tt_k} \times  \frac{\prod\limits_{r=1}^{k+1} (-\gam_{2k,s} - \tt_k \gam_{2k+1,r} + \frac{ic}{2})}{\prod\limits_{r=1}^{k} (\gam_{2k,s} -\tt_k\gam_{2k-1,r} - \frac{ic}{2})} + \\ &
	i \frac{\prod\limits_{r=1}^{k} (-\tt_k\gam_{2k-1,r}) \prod\limits_{r=1}^{k+1} (-\tt_{k}\gam_{2k+1,r})}{(-\tt_{k} \gam_{2k+\delta,j})\prod\limits_{r=1}^{k} (-\frac{ic}{2} - \gam_{2k,r})(-\frac{ic}{2} + \gam_{2k,r})} \Bigg) w_{n} =  i w_{n}
\end{split} \eeqq
The last line is obtained by the use of the following well know identity:
\beq \label{A8}
\sum\limits^{m}_{i=1} \frac{\prod\limits^{m-1}_{j=1} (x_{i} - y_{j})}{\prod\limits_{r \neq i} (x_{i} - x_{r})} = 1
\eeq
where for indeterminates $x_i$ we choose  $2k+2$ variables $\{ \pm \gam_{2k,s}, -\frac{ic}{2}, \tt_k \gam_{2k+\delta} - \frac{ic}{2} \}$, and for indeterminates $y_i$ choose $\{ \tt_k \gam_{2k-1} - \frac{ic}{2}, \tt_k \gam_{2k+1} - \frac{ic}{2}  \}$.

Let us check the equality \rf{A7}. We have
\beqq  \begin{split}
	& J^{\tt_{k+1}}_{k,\delta,j} w_{n} =\\& \bigg(\sum\limits_{s=1}^{k} \frac{\tt_{k+1} ic}{\gam_{2k+\delta,j} - \tt_{k+1}  (\gam_{2k,s}+\frac{ic}{2})} Q_{k,s} + \sum\limits_{s=1}^{k} \frac{\tt_{k+1} ic}{\gam_{2k+\delta,j} + \tt_{k+1}  (\gam_{2k,s}-\frac{ic}{2})} R_{k,s} + \frac{\tt_{k+1} ic}{\gam_{2k + \delta,j}} T_{k} \bigg) w_{n} 
\end{split}\eeqq
\beqq\begin{split}
	& =
	\Bigg( \sum\limits_{s=1}^{k} \frac{\tt_{k+1} ic}{\gam_{2k+\delta,j} - \tt_{k+1} (\gam_{2k,s}+\frac{ic}{2})} \cdot \frac{1}{ic} \cdot \frac{\prod\limits_{r=1}^{k+1} ((\gam_{2k,s}+\frac{ic}{2})^2 - \gam_{2k+1,r}^2)}{2\gam_{2k,s} (\gam_{2k,s}+\frac{ic}{2}) \prod\limits_{r \neq s} (\gam_{2k,s}^2 - \gam_{2k,r}^2)} \times 
\end{split}\eeqq 
\beqq \begin{split} &
	i^{\tt_{k+1}} \times \frac{\prod\limits_{r=1}^{k} (\gam_{2k,s} - \tt_{k} \gam_{2k-1,r} + \frac{ic}{2})}{\prod\limits_{r=1}^{k+1} (-\gam_{2k,s} -\tt_{k} \gam_{2k+1,r} - \frac{ic}{2})} - \\&
	\sum\limits_{s=1}^{k} \frac{\tt_{k+1} ic}{\gam_{2k+\delta,j} + \tt_{k+1} (\gam_{2k,s}-\frac{ic}{2})} \cdot \frac{1}{ic} \cdot \frac{\prod\limits_{r=1}^{k+1} ((\gam_{2k,s}-\frac{ic}{2})^2 - \gam_{2k+1,r}^2)}{2\gam_{2k,s} (\gam_{2k,s}-\frac{ic}{2}) \prod\limits_{r \neq s} (\gam_{2k,s}^2 - \gam_{2k,r}^2)} \times \\& i^{-\tt_{k+1}} \cdot \frac{\prod\limits_{r=1}^{k+1} (-\gam_{2k,s} - \tt_{k} \gam_{2k+1,r} + \frac{ic}{2})}{\prod\limits_{r=1}^{k} (\gam_{2k,s} -\tt_{k} \gam_{2k-1,r} - \frac{ic}{2})} 
	+ i \frac{\prod\limits_{r=1}^{k} (-\tt_{k} \gam_{2k-1,r}) \prod\limits_{r=1}^{k+1} (-\tt_{k} \gam_{2k+1,r})}{(-\tt_{k+1} \gam_{2k+\delta})\prod\limits_{r=1}^{k} (-\frac{ic}{2} - \gam_{2k,r})(-\frac{ic}{2} + \gam_{2k,r})} \Bigg) w_{n} \\ &
	= i \Bigg( 1 - \frac{\prod\limits_{r=1}^{k+1} (\gam_{2k+\delta,j} + \gam_{2k+1,r}) \prod\limits_{r=1}^{k} (\gam_{2k+\delta,j} + \gam_{2k-1,r})}{\gam_{2k+\delta,j} \prod\limits_{r=1}^{k} ((\gam_{2k+\delta,j}+\tt_{k}\frac{ic}{2})^2 - \gam_{2k,r}^2)} \Bigg) w_{n}
\end{split} \eeqq
The last line of the equality is obtained after simplifications of the ratios by  use of the relation $\rf{A8}$, where for inderteminates $x_i$ we substitute  $x$ variables  $\{ \pm \gam_{2k,s}, -\frac{ic}{2}, \tt_{k+1} \gam_{2k+\delta} - \frac{ic}{2} \}$, and for $y_i$ we substitute $\{ \tt_k \gam_{2k-1} - \frac{ic}{2}, \tt_k \gam_{2k+1} - \frac{ic}{2}  \}$.

\bigskip
{\bf Proof of \rf{A2}}. According to Lemmas \ref{lemmaA1} and  \ref{lemmaA2},
for any    $k = 1,\ldots,n$ we have:
 \begin{align} \notag
	& (I_{2k+1,2k} + i I_{2k+1,2k-1})w_{n} = \sum \limits_{\ve = \pm 1} \sum\limits_{j=1}^{k} P^{\ve}_{kj} (1 + i J^{\ve}_{k-1,1,j}) w_{n} = \\
	\label{A10a} &
	\sum\limits_{j=1}^{k} P^{\tt_{k}}_{kj} \frac{\prod\limits_{r=1}^{k} (\gam_{2k-1,j} + \gam_{2k-1,r}) \prod\limits_{r=1}^{k-1} (\gam_{2k-1,j} + \gam_{2k-3,r})}{\gam_{2k-1,j} \prod\limits_{r=1}^{k-1} ((\gam_{2k-1,j}+\tt_{k-1}\frac{ic}{2})^2 - \gam_{2k-2,r}^2)} w_{n} = \end{align}\beqq\begin{split} &
	- \frac{1}{ic} \sum\limits_{j=1}^{k}  \frac{\prod\limits_{r=1}^{k-1} (\gam_{2k-1,j} + \tt_{k} (\gam_{2k-2,r} + \frac{ic}{2})) \prod\limits_{r=1}^{k} (\gam_{2k-1,j} - \tt_{k} (\gam_{2k,r} - \frac{ic}{2}))}{2\prod\limits_{r \neq j} (\gam_{2k-1,j} - \gam_{2k-1,r})(\gam_{2k-1,j} + \gam_{2k-1,r} + \tt_{k} ic)} e^{\tt_{k} ic \partial_{\gam_{2k-1,j}}} \times \\ &
	2 \frac{\prod\limits_{r \neq j}^{k} (\gam_{2k-1,j} + \gam_{2k-1,r}) \prod\limits_{r=1}^{k-1} (\gam_{2k-1,j} + \gam_{2k-3,r})}{\prod\limits_{r=1}^{k-1} ((\gam_{2k-1,j}+\tt_{k-1}\frac{ic}{2})^2 - \gam_{2k-2,r}^2)} w_{n} =\\ &
	- \frac{1}{ic} \sum\limits_{j=1}^{k} \frac{\prod\limits_{r=1}^{k-1} (\gam_{2k-1,j} + \gam_{2k-3,r} + \tt_{k}ic)}{\prod\limits_{r \neq j} (\gam_{2k-1,j} - \gam_{2k-1,r})} \cdot \frac{\prod\limits_{r=1}^{k} \left(\gam_{2k-1,j} + \tt_{k+1} \gam_{2k,r} + \tt_{k} \frac{ic}{2}\right)}{\prod\limits_{r=1}^{k-1} \left(\gam_{2k-1,j} + \tt_{k+1} \gam_{2k-2,r} + \tt_{k} \frac{ic}{2}\right)} e^{\tt_{k} ic \partial_{\gam_{2k-1,j}}} w_{n} 
\end{split} \eeqq
In the product, presenting the function $w_{n}$, only the factors  $W_{k-1}^{\tt_{k}}$ and  $W_{k}^{\tt_{k+1}}$  depend on the variables $\gam_{2k-1,j}$.
Thus, using \rf{w2}, functional relations on the Euler Gamma function and \rf{A8}, we get
\beq \label{A10} \begin{split}
	&(I_{2k+1,2k} + i I_{2k+1,2k-1})w_{n} =  \frac{\tt_{k+1}}{ic} \sum\limits_{j=1}^{k} \frac{\prod\limits_{r=1}^{k-1} (\gam_{2k-1,j} + \gam_{2k-3,r} + \tt_{k}ic)}{\prod\limits_{r \neq j} (\gam_{2k-1,j} - \gam_{2k-1,r})} w_{n} = \frac{\tt_{k+1}}{ic} w_{n}
\end{split} \eeq
\hfill{$\Box$}

{\bf Proof of \rf{A2a}}.
Again,  according to Lemmas \ref{lemmaA1} and  \ref{lemmaA2},
for any    $k = 1,\ldots,n$ we have:        
\beq \begin{split}
	& (-I_{2k+2,2k-1} + i I_{2k+2,2k})w_{n} = \sum \limits_{\ve = \pm 1} \sum\limits_{j=1}^{k} P^{\ve}_{kj} J^{\ve}_{k,-1,j} (J^{\ve}_{k-1,1,j} - i) w_{n} = \\ &
	-i \sum\limits_{j=1}^{k} P^{\tt_{k}}_{kj} J^{\tt_{k}}_{k,-1,j} \frac{\prod\limits_{r=1}^{k} (\gam_{2k-1,j} + \gam_{2k-1,r}) \prod\limits_{r=1}^{k-1} (\gam_{2k-1,j} + \gam_{2k-3,r})}{\gam_{2k-1,j} \prod\limits_{r=1}^{k-1} ((\gam_{2k-1,j}+\tt_{k-1}\frac{ic}{2})^2 - \gam_{2k-2,r}^2)} w_{n}. 
\end{split} \eeq
Since the operator $J^{\tt_{k}}_{k,-1,j}$ contains shifts only of the variables $\gam_{2k,s}$, we can rewright the result as 
\beq \label{A11}
-i \sum\limits_{j=1}^{k} P^{\tt_{k}}_{kj} \frac{\prod\limits_{r=1}^{k} (\gam_{2k-1,j} + \gam_{2k-1,r}) \prod\limits_{r=1}^{k-1} (\gam_{2k-1,j} + \gam_{2k-3,r})}{\gam_{2k-1,j} \prod\limits_{r=1}^{k-1} ((\gam_{2k-1,j}+\tt_{k-1}\frac{ic}{2})^2 - \gam_{2k-2,r}^2)} J^{\tt_{k}}_{k,-1,j} w_{n} \eeq
Now  Lemma \ref{lemmaA2} says that 
$$J^{\tt_{k}}_{k,-1,j}w_n=iw_n.$$ 
Thus $(-I_{2k+2,2k-1} + i I_{2k+2,2k})w_{n}$ equals to 
\beq\label{A13}
= \sum\limits_{j=1}^{k} P^{\tt_{k}}_{kj} \frac{\prod\limits_{r=1}^{k} (\gam_{2k-1,j} + \gam_{2k-1,r}) \prod\limits_{r=1}^{k-1} (\gam_{2k-1,j} + \gam_{2k-3,r})}{\gam_{2k-1,j} \prod\limits_{r=1}^{k-1} ((\gam_{2k-1,j}+\tt_{k-1}\frac{ic}{2})^2 - \gam_{2k-2,r}^2)} w_{n} =  \frac{\tt_{k+1}}{ic} w_{n}
\eeq
 The latter equality was proved during the derivation of \rf{A10} from \rf{A10a}.

\hfill{$\Box$}

This ends the proof of Theorem \ref{theorem1}. \hfill{$\Box$}
\medskip

{\bf Remark}. Analyzing the proof of Theorem \ref{theorem1} we see that the derivations over variables $\gam_{2k-1,j}$ enter the game only in the last stage of calculations. Moreover, we can freely add to Whittaker vectors factors of the form
\beqq e^{i\frac{\alpha_j}{c} \delta_{2j-1}},\qquad\text{where}\qquad \delta_{2j-1}=\sum_i\gam_{2j-1,i}=(-1)^{j+1}\sum_i\g_{2j-1,i}\eeqq
where $\alpha $ is arbitrary real number. This does not affect to the convergence of integrals and does not change the action of Cartan subalgebra.
 In Toda equation we earn thus arbitrary positive constants $c_j=e^{\alpha_j}$ at exponentials $e^{x_{j-1}-x_j}$ which can be equivalently obtained by successice shifts of the variables $x_j$.
 	 
\renewcommand*{\thesection}{B}
\setcounter{equation}{0}
\section{Action of Cartan subalgebra}
It is sufficient to calculate $F_{k,k}w_{n}$, and then use the automorphism $\tau$.
\beq \begin{split}
	& -iI_{2k,2k-1} w_{n} =\\& - \Bigg( \frac{1}{c} \sum_{j=1}^{k-1} \frac{\prod\limits_{r=1}^{k} ((\gam_{2k-2,j}+\frac{ic}{2})^2 - \gam_{2k-1,r}^2)}{2\gam_{2k-2,j} (\gam_{2k-2,j}+\frac{ic}{2}) \prod\limits_{r \neq j} (\gam_{2k-2,j}^2 - \gam_{2k-2,r}^2)}  i^{\tt_{k-1}} \frac{\prod\limits_{r=1}^{k-1} (\gam_{2k-2,j} - \tt_{k-1} \gam_{2k-3,r} + \frac{ic}{2})}{\prod\limits_{r=1}^{k} (-\gam_{2k-2,j} - \tt_{k-1} \gam_{2k-1,r} - \frac{ic}{2})}  \\ & + \frac{1}{c} \sum_{j=1}^{k-1} \frac{\prod\limits_{r=1}^{k-1} ((\gam_{2k-2,j}-\frac{ic}{2})^2 - \gam_{2k-3,r}^2)}{2\gam_{2k-2,j} (\gam_{2k-2,j}-\frac{ic}{2}) \prod\limits_{r \neq j} (\gam_{2k-2,j}^2 - \gam_{2k-2,r}^2)} i^{- \tt_{k-1}} \frac{\prod\limits_{r=1}^{k} (- \gam_{2k-2,j} - \tt_{k-1} \gam_{2k-1,r} + \frac{ic}{2})}{\prod\limits_{r=1}^{k-1} (\gam_{2k-2,j} - \tt_{k-1} \gam_{2k-3,r} - \frac{ic}{2})} - \\ & 
	- \frac{1}{ic} \frac{\prod\limits_{r=1}^{k-1} \gam_{2k-3,r} \prod\limits_{r=1}^{k} \gam_{2k-1,r}}{\prod\limits_{r=1}^{k-1} (\gam_{2k-2,r}^2 + \frac{c^2}{4})} \Bigg) w_{n} = \end{split}\eeq
\beq\begin{split} 
	= - &\frac{1}{ic} \Bigg( \sum_{j=1}^{k-1} \frac{\prod\limits_{r=1}^{k} (\gam_{2k-2,j} - \tt_{k-1} \gam_{2k-1,r}+\frac{ic}{2}) \prod\limits_{r=1}^{k-1} (\gam_{2k-2,j} - \tt_{k-1} \gam_{2k-3,r} + \frac{ic}{2})}{2\gam_{2k-2,j} (\gam_{2k-2,j}+\frac{ic}{2}) \prod\limits_{r \neq j} (\gam_{2k-2,j}^2 - \gam_{2k-2,r}^2)} + \\ &  \sum_{j=1}^{k-1} \frac{\prod\limits_{r=1}^{k-1} (-\gam_{2k-2,j} - \tt_{k-1} \gam_{2k-3,r} + \frac{ic}{2}) \prod\limits_{r=1}^{k} (-\gam_{2k-2,j} - \tt_{k-1} \gam_{2k-1,r} + \frac{ic}{2})}{2\gam_{2k-2,j} (\gam_{2k-2,j}-\frac{ic}{2}) \prod\limits_{r \neq j} (\gam_{2k-2,j}^2 - \gam_{2k-2,r}^2)} + \\ &
	 \frac{\prod\limits_{r=1}^{k-1} (-\tt_{k-1} \gam_{2k-3,r}) \prod\limits_{r=1}^{k} ( - \tt_{k-1} \gam_{2k-1,r})}{\prod\limits_{r=1}^{k-1} (-\frac{ic}{2} - \gam_{2k-2,r}) (-\frac{ic}{2} + \gam_{2k-2,r})} \Bigg) w_{n} =\\& \frac{\tt_{k-1}}{ic} \left( \sum\limits^{k}_{j=1} \gam_{2k-1,j} + \sum\limits^{k-1}_{j=1} \gam_{2k-3,j} - \tt_{k-1} (k-1)ic \right).
\end{split} \eeq
Here we again use the identity \rf{A8}
where for $x_i$ we substitute  $2k-1$ variables $\{ \pm \gam_{2k-2,s}, -\frac{ic}{2} \}$, а and for $y_i$ use $\{ \tt_{k-1} \gam_{2k-1} - \frac{ic}{2}, \tt_{k-1} \gam_{2k-3} - \frac{ic}{2}  \}$.
\hfill{$\Box$}

\end{document}